\documentclass[11pt,a4paper]{article} 

\usepackage{amssymb}

\def\ccp{\mathcal C}
\def\aap{\mathcal A}

\def\aho{\alpha}
\def\cho{\kappa}
\def\rhom{\varrho}
\def\fraka{D}
\def\mao{\Theta}
\def\gga{m}
\def\Gga{g}

\newcommand{\sa}{\,|\,}
\newcommand{\cf}{\mathcal{F}}
\newcommand{\cg}{\mathcal{G}}
\newcommand{\E}{\mathbb{E}}
\newcommand{\dd}{\mathrm{d}}
\newcommand{\ee}{\mathrm{e}}
\renewcommand{\P}{\mathbb{P}}

\newcommand{\R}{\mathbb{R}}

\renewcommand{\le}{\leqslant}
\renewcommand{\ge}{\geqslant}

\newcommand{\vide}{\mbox{\O}}
\newcommand{\un}{\mathbf{1}}

\newtheorem{theorem}{Theorem}

\newtheorem{lemma}{Lemma}
\newtheorem{corollary}[lemma]{Corollary}

\newtheorem{proposition}[lemma]{Proposition}

\def\proof#1{\noindent\trivlist\item[\hskip\labelsep{\bf #1}]\ignorespaces}
\def\endproof{\hfill $\square\qquad$\endtrivlist}

\newtheorem{remark}[lemma]{Remark}
\newenvironment{rbox}{\begin{remark}\begin{em}}{\hfill$\square\qquad$\end{em}\end{remark}}

\begin{document}
\title{Invariance principle for the coverage rate of genomic physical mappings}
\author{Didier Piau}
\date{\em Universit\'e Lyon 1}
\maketitle

\thispagestyle{empty}

\begin{abstract}
We study some stochastic models of physical mapping of genomic
sequences. Our starting point is a global construction of the process
of the clones and of the process of the anchors which are used to map
the sequence. This yields explicit formulas for the moments of the
proportion occupied by the anchored clones, even in inhomogeneous
models. This also allows to compare, in this respect, inhomogeneous models
to homogeneous ones. Finally, for homogeneous models, we provide
nonasymptotic bounds of the variance and we prove functional
invariance results.
\end{abstract}

\bigskip

Date: \today.

Abbreviated title:
Invariance of physical mappings




\bigskip

\section*{Introduction}

The goal of the projects of genomic physical mapping is to reconstruct
almost completely the sequence of a genome, starting from a multitude of
exactly sequenced fragments, which are called clones. One approach to
the reconstruction of the overall positions of these clones in the
complete genomic sequence uses so-called anchors. These are short,
exactly sequenced, portions of the genome which are assumed to appear
only once in the full genomic sequence.  An anchored clone is a clone
which contains an anchor. In this paper, we assume that the positions
of the anchors, hence of the anchored clones, are exactly known.
Maximal connected unions of anchored clones are called islands or,
more exactly, anchored islands, aka contigs. The complement of the
islands is called the ocean.  When suitably rescaled, the full genomic
sequence is identified with (a portion of) the real line, the anchors
are identified with points, and the clones and the islands are
identified with intervals.

The overall quality of the reconstitution of a given genomic sequence
depends obviously on the number of islands, on their length, and on the
proportion of the sequence which is occupied by the ocean, among other
characteristics of the project. One hopes that the islands are as few
and as long as possible, and that the proportion occupied by the ocean
is as low as possible.  Arratia et al.\ (1991) introduced a stochastic
model of physical mapping, where the positions of the right ends of
the clones and the positions of the anchors are distributed according
to independent
homogeneous Poisson processes on the real line, and where the lengths of the
clones are random, i.i.d., and independent of everything else.  For this
model, Arratia et al.\ computed the mean values of the three
quantities of interest that we mentioned above.
For related studies, see Lander and
Waterman\ (1988), Ewens et al.\ (1991), and Grigoriev\ (1993).

Motivated by the fact that actual genomic sequences do not fulfill the
homogeneity hypotheses which underly the stochastic model introduced 
by Arratia et
al., Schbath (1997) and Schbath et al.\ (2000)
extend this setting in two directions.  In both papers, the
independence properties of the model remain, but Schbath (1997) studies
the case when the
intensities of the Poisson processes which generate the positions of
the clones and the positions of the anchors may depend on their
respective positions along the genome, and Schbath et al.\ (2000) study
the case when the distributions of the
lengths of the clones may depend on their respective positions along
the genome.  
In these two wider contexts, these papers
provide expressions of the mean value of the number of islands, of
the mean value of the proportion occupied by the ocean, and, under an
additional technical hypothesis, of the mean value of the length of
the islands.

In the present paper, we pursue the study of this class of models. As a
first contribution, we consider the class of models where the
Poisson process of the clones, the Poisson process of the anchors, and
the distributions of the lengths of the clones can all be
inhomogeneous simultaneously. To give a flavor of our results in this direction, we
state proposition~\ref{p.1} below, which extends formulas of the
papers mentioned above, for the mean value of the number of
clones and for the mean value of the number of anchored clones which
cover a point.

To state proposition~\ref{p.1}, we introduce the measure $c(\dd x)$ on the real line as the intensity
measure of the Poisson process of the (right ends of the) clones, the
measure $a(\dd x)$ on the real line as the intensity measure of the
Poisson process of the anchors, and, for every $x$ on the real line,
the random variable 
$L_x$ as the length of a clone whose right end is at position $x$, and
we refer to section~\ref{s.gm} for more precise definitions of these objects.

For every $x$ on the real line, $n_\ccp(x)$ denotes the number of
clones which contain the point $x$, and $n_\aap(x)$ denotes the number
of anchored clones which contain the point $x$.

\begin{proposition}[General case]
\label{p.1}
\textbf{\em (1)}
The random variable $n_\ccp(x)$ follows the Poisson
distribution whose mean value is given by the expression
$$
\E(n_\ccp(x))=\int_x^{+\infty}c(\dd z)\,\P(L_z\ge z-x).
$$
\textbf{\em (2)}
The mean value of $n_\aap(x)$ is  given by the expression
$$
\E(n_\aap(x))=\int_x^{+\infty}c(\dd
z)\int_{z-x}^{+\infty}\P(L_z\in\dd t)\,\left(1-\ee^{-a([z-t,z])}\right).
$$ 
Here $a([z-t,z])$ denotes the mesure of the interval $[z-t,z]$ with
respect to the measure $a(\dd x)$.

\textbf{\em (3)}
The distribution of the random variable $n_\aap(x)$ is not Poisson. 
More specifically, either $n_\aap(x)=0$ almost surely, or 
the variance of $n_\aap(x)$ is
strictly greater than its mean value.  
\end{proposition}

In
actual physical mapping projects, the condition that $n_\aap(x)=0$ almost surely is never fulfilled. On the mathematical side, this would
correspond to degeneracies such as the fact that $L_z\le z-x$ almost
surely for every $z\ge x$ and/or the fact that the intensity of the
anchors is zero on a suitable neighborhood of $x$.

The homogeneous case is when  $c(\dd x)=\cho\,\dd x$ and
$a(\dd x)=\aho\,\dd x$ for two given positive constants $\cho$ and
$\aho$, 
and when every $L_x$ is distributed like a given random
variable $L$.
The specialization of proposition~\ref{p.1} to the homogeneous case is as follows.

\begin{corollary}[Homogeneous case]
\label{c.h}
In the homogeneous case with parameters $\cho$, $\aho$, and $L$, the
mean values of $n_\ccp(x)$ and $n_\aap(x)$ do not depend on the point $x$
and
are given by the expressions
$$
\E(n_\ccp)=\cho\,\E(L),
\quad
\E(n_\aap)=\cho\,\E(L\,(1-\ee^{-\aho L})).
$$
\end{corollary}

More importantly than the slight generalizations above,
our second contribution is to provide explicit formulas for 
the higher moments of these quantities in the general model with
variable intensities. In the homogeneous case, our results imply,
for instance, that the proportion of a large genomic sequence occupied
by the ocean is
asymptotically Gaussian, see theorem~\ref{t.a} below.

\begin{theorem}[Homogeneous case]
\label{t.a}
Consider the homogeneous case with parameters $\cho$, $\aho$, and $L$, and assume that 
$L$ is square integrable.  For any
positive $G$, let the random variable $O_G$ denote the measure of the intersection of the
ocean with any interval of length $G$, for instance the interval $[0,G]$, and let $\sigma^2(O_G)$ denote the
variance of $O_G$.

\textbf{\em (1)}
There exists a positive constant $\rhom<1$ such that
$\E(O_G)=\rhom\,G$ for every nonnegative $G$.

\textbf{\em (2)}
There exists finite positive constants $\nu$ and $\lambda$ such that,
for every nonnegative $G$,
$$
\nu\,G-\lambda\le\sigma^2(O_G)\le\nu\,G.
$$ 
Hence
$\sigma^2(O_G)\sim\nu\,G$ when $G\to\infty$.
Furthermore, the function $G\mapsto\sigma^2(O_G)$ is convex, and
$\sigma^2(O_G)-\nu\,G+\lambda\to0$  when $G\to\infty$.

\textbf{\em (3)}
For every positive $G$, let
$\mao_G$ denote the random 
process, indexed by the real numbers $0\le t\le1$, and defined by
$$
\mao_G(t):=(O_{Gt}-\rhom\,Gt)/\sqrt{\nu\,G}.
$$
When $G\to\infty$, the process $\mao_G$  
converges in distribution to a standard 
Wiener process on the space of continuous functions on $[0, 1]$, 
equipped with the metric of the
uniform convergence.

\textbf{\em (4)}
The constants $\rhom$, $\nu$ and $\lambda$ above
can be written explicitly as integrals which
involve the parameters $\cho$, $\aho$, and the distribution of $L$.
\end{theorem}

The starting point of our results is a global construction of the
clones, the anchors, and the islands, using a single Poisson
process. We expose this global construction in section~\ref{s.gm}.  We
provide alternative descriptions of this process, locating for
instance the clones by their left ends instead of their right ends.  A
natural conjecture in this setting is that the homogeneous model would
be the only one invariant by the symmetry of the real line, but we
disprove this.  In section \ref{s.m}, we rewrite in our general
setting various formulas due to Arratia et al.\ or to Schbath or to
Schbath et al.
Section \ref{s.hm} provides explicit formulas for every moment of the
proportion of the real line which is occupied by the ocean in the general case and provides
rather sharp bounds of the variance in the homogeneous case. Finally,
section~\ref{s.ihc} proves the invariance result stated in theorem~\ref{t.a}
above, in the homogeneous case.  On our way, we provide asymptotics of the
variance when the number of clones is vanishingly small and we build
comparison tools that yield effective upper and lower bounds in some
inhomogeneous cases.

\paragraph{Acknowledgements}
We wish to thank Sophie Schbath for an introduction to this
subject and for instructive discussions, Julien Michel for the
positive association argument used in section~\ref{ss.pd}, and an
anonymous referee for a careful reading of the first version of this paper.

\section{Global model}
\label{s.gm}

In this section, we build the clones, the anchors and the islands from
a single Poisson process. Sections~\ref{ss.ad} and \ref{ss.ns}  are
not used in the rest of the paper and may be omitted on a first
reading.

\subsection{Clones}
\label{ss.c}
Let $\R$ denote the real line and $\R^+:=[0,+\infty)$  the
  nonnegative half line.
Let $c(\dd x)$ denote the intensity measure of the Poisson process of
the right ends of the clones. Assume that, when the 
right end of a clone is located at $x$, its length follows
the distribution of a given random variable $L_x$.  We represent the
clone which covers exactly the interval $[x-t,x]$ of length $t\ge0$ by
the point $(x,t)$ in $\R\times\R^+$.  The distribution of the clones is described by a 
Poisson process $\ccp$ on $\R\times\R^+$ of
intensity measure $\gga$, with
$$
\gga(\dd x\,\dd t):=c(\dd x)\,\P(L_x\in\dd t).
$$ 
In other words, $\ccp$ is a random subset of $\R\times\R^{+}$, which is
almost surely locally finite, 
and such that the following holds.
For every Borel subsets $D$ and $D'$ of 
$\R\times\R^{+}$ such that $D\cap D'$ is empty, 
the random number of points of $\ccp$ in $D$ and the random number of $\ccp$ in $D'$ 
are independent. Furthermore, for every Borel subset $D$ of 
$\R\times\R^{+}$, the number of points of $\ccp$
in $D$ is a 
Poisson random variable of mean value $\gga(D)$.

In fact, the intensity measure $\gga$ can be any Borel measure on $\R\times\R^+$ with a
locally finite first marginal $c$, given by
$$
c(\dd x):=\int_{t\ge0}\gga(\dd x\,\dd t).
$$
That is, one assumes that $c([-G,G])=\gga([-G,G]\times\R^+)$ is
finite for every finite positive $G$.
The assumption that $c$ is locally finite ensures that the
distribution 
of $L_x$ is well
defined, and given by the Radon-Nikodym derivative
$$
\P(L_x\in\dd t):=\gga(\dd x\,\dd t)/c(\dd x).
$$

\subsection{Alternative descriptions of the clones}
\label{ss.ad}
At first sight, it may seem rather arbitrary to locate the position of
a clone by its right endpoint, rather than by its midpoint or by its
left endpoint.  In fact, these alternative descriptions are also
characterized by Poisson processes, albeit possibly with different
intensities.  For instance, using the couple $(y,t)$ to describe the
clone $[y,y+t]$ yields a Poisson process on $\R\times\R^+$ of
intensity measure $\gga'$, with
$$
\gga'(\dd y\,\dd t):=:c'(\dd y)\,\P(L'_y\in\dd t).
$$ 
One obtains  $\gga'$ from $\gga$, or rather, one obtains 
$c'(\dd y)$ and the distributions of the random variables $L'_y$ 
from $c(\dd x)$ and
from the distributions of the random variables $L_x$, as follows.
For any nonnegative test function $\Phi$, the expected value of the sum over
every clone $[y,x]$ of $\Phi(y,x)$ reads
$$
\E\left(\sum_{[y,x]\,\mbox{\small clone}}\Phi(y,x)\right)=
\int\!\!\!\int\gga(\dd x\,\dd t)\,\Phi(x-t,x)=
\int\!\!\!\int\gga'(\dd y\,\dd t)\,\Phi(y,y+t).
$$
In other words, one asks that
$$
\int c(\dd x)\,\E(\Phi(x-L_x,x))=
\int c'(\dd y)\,\E(\Phi(y,y+L'_y)).
$$
Since this equality holds for every test function $\Phi$, 
this implies that $c'(\dd y)$ and the distributions of the random
variables $L'_y$ are given by
\begin{eqnarray*}
c'(\dd y)
& = &
\int_{y}^{+\infty}c(\dd x)\,\P(L_x\in x-\dd y),
\\
\P(L'_y\in\dd t)
& = &
\P(L_{y+t}\in\dd t)\,c(t+\dd y)/c'(\dd y).
\end{eqnarray*}
Similar formulas give the intensity measures associated to the
description of a clone by its midpoint and by its length, or by its
two endpoints.

\subsection{On the (non)specificity of the homogeneous clones}
\label{ss.ns}
Based upon the preceding section, the reader might be led to
believe that the homogeneous model is privileged with respect to the
transformations of the intensity measure $\gga(\dd x\,\dd t)$ into
the intensity measure $\gga'(\dd y\,\dd t)$ and of $\gga'(\dd
y\,\dd t)$ into $\gga(\dd x\,\dd t)$.  To wit, if the
intensity $c(\dd x)$ and the distributions of the random variables
$L_x$ are invariant by the translations of the real line, so are the
intensity $c'(\dd y)$ and the distributions of the random variables
$L'_y$. Thus, in the homogeneous case, $c(\dd x)=c'(\dd x)=\cho\,\dd
x$ and the distributions of every $L_x$ and every $L'_y$ do coincide.

Our goal in this section is to point out that there are other cases
where the two intensity measures $\gga$ and $\gga'$ coincide.
To build such examples, we need to introduce, for every $x$ on the real line, 
the unit interval $U_x$ which is centered
at $x$, that is,
$$
U_x:=[x-1/2,x+1/2).
$$
Let
$B_0$ denote the union of the intervals $U_{2k}$ 
for every integer $k$, and let $B_1$ denote
its complement.
Let  
$u_0(\dd x)$ denote a finite measure on $U_0$, and
    $u_1(\dd x)$ a finite measure on $U_1$. Let $c(\dd x)$ denote the unique
measure on $\R$ which is invariant by the translation $x\mapsto x+2$
and whose restrictions to $U_0$ and to $U_1$ are $u_0(\dd x)$ and
$u_1(\dd x)$, respectively.
Thus, $c=c_0+c_1$ with, for $i=0$ and for $i=1$,
$$
c_i(\dd x):=\sum_ku_i(2k+\dd x),
$$
where both sums run over every integer $k$.
In other words, $c(\dd x)$ can be any locally finite measure on $\R$,
invariant by the translation $x\mapsto x+2$, and the measure $c_0(\dd x)$,
respectively the measure $c_1(\dd x)$, denotes the restriction of $c(\dd x)$ to
$B_0$, respectively the restriction of $c(\dd x)$ to $B_1$.

Assume finally that $L_x=2$ with full probability when $x$ is in
$B_0$, and that $L_x=4$ with full probability when $x$ is in $B_1$.
Since $L_x$ is always an even integer, the endpoints of a given clone are
either both in $B_0$ or both in $B_1$.  Using this
remark, one can check that $\gga=\gga'$. Besides, the process which locates the
clones by their midpoint is given by a similar intensity measure,
choosing with full probability the length $4$ when the midpoint
belongs to $B_0$, and choosing with full probability the length $2$
when the midpoint belongs to $B_1$.

In the example above, the distributions of the lengths are discrete,
hence the measure $\gga(\dd x\,\dd t)$ is singular with respect to the Lebesgue
measure.
However, the same idea can be adapted to produce examples where $\gga(\dd x\,\dd t)$ is
absolutely continuous.
To see this, introduce the Poisson process which describes a clone
$[y,x]$ by its endpoints $(y,x)$, and assume that the intensity measure
$\gga^*$ of this Poisson process is
$$
\gga^*(\dd y\,\dd x):=\dd y\,\dd x\sum_{k}\un\{(y,x)\in U_{2k}\times U_{2k+2}\}+
\un\{(y,x)\in U_{2k-1}\times U_{2k+3}\},
$$ 
where the sum runs over every integer $k$.
In words, the left endpoints and the right endpoints of the clones both have homogeneous
intensity measures, and both endpoints of a clone belong to $B_0$ or both endpoints belong to
$B_1$. Furthermore, given that the left endpoint $y$ belongs to $B_0$, the
right endpoint $x$ is uniformly distributed over the next unit interval of $B_0$
to the right of $y$, that is, over the connected component of $B_0$ which
contains $y+2$. Given that the left endpoint $y$ belongs to $B_1$, the
right endpoint $x$ is uniformly distributed over the second next 
unit interval of $B_1$
to the right of $y$, that is, over the connected component of $B_1$
which contains $y+4$.

In this new example, the measure $\gga(\dd x\,\dd t)$ is as follows. 
The intensity $c(\dd x)$ is the Lebesgue
measure. The length $L_x$ is uniformly distributed over $U_{x-2k}$ when
$x$ is in $U_{2k+2}$, and $L_x$ is uniformly distributed over
$U_{x-2k+1}$ when $x$ is in $U_{2k+3}$.  The support of the
distribution of $L_x$ is a unit subinterval of the interval $[1,3]$
when $x$ is in $B_0$, and it is a unit subinterval of the interval
$[3,5]$ when $x$ is in $B_1$, hence the distribution of $L_x$ cannot
be the same for every $x$ on the real line.  Finally, $\gga=\gga'$ because $\gga^*$
is invariant by the symmetries of the real line, since these exchange the left
endpoints and the right endpoints of the clones while leaving their lengths
unchanged.

\subsection{Anchors}
\label{ss.a}
In this section and in the rest of the paper, we come back to the
$(x,t)$ Poisson process of intensity $\gga$, which represents the
clones by their right endpoint and by their length.

The anchors are described by a Poisson process $\aap$ on the real line, with intensity
$a(\dd x)$, independent of the Poisson process $\ccp$ of
the clones which we defined in section~\ref{ss.c}. Thus, for every Borel subset $D$ of the real line, the number of
anchors in $D$ is a random variable whose distribution is Poisson with
mean value $a(D)$, and the number of anchors in the Borel sets $D$ and
$D'$ are independent random variables as soon as $D\cap D'$ is empty.

For every subset
$\fraka$ of the real line, let $I(\fraka)$ denote the cone of
influence of $\fraka$ in $\R\times\R^+$.
This is the set of clones $(x,t)$ which become anchored clones when
every point of $\fraka$ becomes an anchor. Thus,
$$
I(\fraka):=\{(x,t)\in\R\times\R^+\,;\,[x-t,x]\cap\fraka\neq\vide\}.
$$
For every measurable $\fraka$, the process $\ccp_{\fraka}:=\ccp\cap I(\fraka)$ of the clones that are
anchored by 
 $\fraka$ is deduced from $\ccp$ 
by erasing some clones, hence each $\ccp_{\fraka}$ is indeed a Poisson
process whose intensity measure $\gga_D$ on $\R\times\R^+$ is the restriction
of the original intensity measure $\gga$ to the set  $I(\fraka)$,
that is,
$$
\gga_{\fraka}(\dd x\,\dd
t):=\un\{(x,t)\in I(\fraka)\}\,\gga(\dd x\,\dd t). 
$$
For every locally finite subset $D$ of the real line, let
$\P^{\fraka}$ denote the conditioning of $\P$ by the event
$\{\ccp=\fraka\}$. Finally, let $\ccp_\aap$ denote the process of the anchored
clones, that is
$$
\ccp_\aap:=\{(x,t)\in \ccp\,;\,[x-t,x]\cap\aap\neq\vide\}=\ccp\cap I(\aap).
$$

\subsection{Clones+anchors}
\label{ss.ca}
One can, and we shall, simultaneously generate the processes $\ccp$, $\aap$ and $\ccp_\aap$
from a unique Poisson process, as follows.
Let  $M:=\R^+\cup\{*\}$, where $*$ denotes any point which is not in $\R^+$.
We endow 
the set $M$ with the smallest $\sigma$-algebra 
which contains the
Borel sets of $\R^+$ and the singleton $\{*\}$. We endow 
the set $\R\times M$ with the product $\sigma$-algebra of the Borel
$\sigma$-algebra of $\R$ and of this $\sigma$-algebra of $M$.
Finally, we introduce a Poisson process on
$\R\times M$ with intensity
$$
\Gga(\dd x\,\dd t):=\gga(\dd x\,\dd t)+a(\dd x)\,\delta_{*}(\dd t).
$$
We call this Poisson process the global process.
The point $(x,t)$ with $t$ in $\R^+$ represents the clone
$[x-t,x]$ and the point $(x,*)$ represents the anchor at $x$. 
The restriction of the global process to the domain $\R\times\R^+$
yields the process of the clones described in section~\ref{ss.c},
since its intensity, which is the
restriction
of
$\Gga$ to $\R\times\R^+$, is
$\gga(\dd x\,\dd t).$
Likewise, the projection $(x,*)\mapsto x$
on the real coordinate of the restriction of
the global process to the domain
$\R\times\{*\}$ yields the
process of the anchors described in section~\ref{ss.a},
since its intensity is $a(\dd x)$.
Finally, the process of the clones and the process of the anchors 
are indeed independent
since they are realized as the restrictions of the global Poisson process to the domains 
$\R\times\R^+$ and $\R\times\{*\}$, which are disjoint subsets of $\R\times M$.

Proposition~\ref{p.t} below and proposition~\ref{p.1} and
corollary~\ref{c.h}
in our introduction
 follow from the construction above. The proofs are simple adaptations
 of the proofs given by Arratia et al., Schbath, and Schbath et al., 
hence we omit them.

\begin{proposition}[General case]
\label{p.t}
With respect to $\P$, $\aap$ and $\ccp$ are independent Poisson processes.
For every locally finite $\fraka$, with respect to $\P^{\fraka}$, $\ccp_{\fraka}$ is a
Poisson process. With respect to $\P$, $\ccp_\aap$ is not a Poisson process.
\end{proposition}

\subsection{Ocean}
\label{ss.o}
Recall that the ocean $O$ is the complement of the union of the anchored
islands. 
For every Borel set $D$ of the real line, let $O(D)$ denote the measure of $O\cap
D$.
For every positive real number $G$, let $O_G:=O([0,G])$. 
For every Borel set $Z$ of the real line, let $$r(Z):=\P(Z\subset O).$$ 
For every $n\ge1$ and every real numbers $z_1$, \ldots, $z_n$, let
$$
r(z_1,\ldots,z_n):=r(\{z_1,\ldots,z_n\})=\P(z_{1}\in O,\ldots,z_{n}\in O).
$$
For instance, $r(z)$ is the probability that $z$ belongs to no
anchored clone. Hence $r(z)$ may depend on $z$ but $r(z)$ corresponds
to $r(0)$ if the process of the clones and the process of the anchors are both shifted by $z$.
Lemma~\ref{l.t} below stems from the definitions.

\begin{lemma}
\label{l.t}
For every Borel set $D$ of the real line and every integer $n\ge1$,
$$
\E(O(D)^n)=\int_{D^n}r(z_1,\ldots,z_n)\,\dd z_1\ldots\dd z_n.
$$
For instance,
$$
\E(O_G)=\int_{0}^{G}r(z)\,\dd z,
\quad
\E(O_G^2)=\int_{0}^{G}\!\!\int_{0}^{G}r(z,z')\,\dd z\,\dd z'.
$$
\end{lemma}


\section{First moments}
\label{s.m}
This section is mainly a rephrasing of results of Arratia et al.\ and Schbath
et al. Our only contribution here is to include both inhomogeneities
simultaneously in the
results, namely, the inhomogeneities of the
 lengths of the clones on the one hand, 
and the inhomogeneities of the positions of the right ends of the clones and of
the anchors on the other hand.
We are interested in $r(z)$, which describes locally the mean value of the
proportion of the real line which is occupied by
the ocean.

\begin{lemma}
Let $J(x,y)$ denote the probability of the event that two
points $x$ and $y$ such that $x\le y$
belong to no common clone.
Then
$$
J(x,y):=\exp\left(-\int_y^{+\infty}\P(L_t\ge t-x)\,c(\dd t)\right).
$$
\end{lemma}

Caution: we renamed
$J(x,x+t)$ the expression $J(x,t)$ of the papers mentioned above.

\begin{lemma}
For every $z$, the joint law of the positions $x$ and $y$ of the anchors which
are the  closest of $z$ to the left and to the right, respectively, is
$A^-(z,\dd
x)\,A^+(z,\dd y)$,
where
$$
A^-(z,\dd x):=A(x,z)\,a(\dd x)
\quad\mbox{et}\quad
A^+(z,\dd y):=A(z,y)\,a(\dd y).
$$ 
For every points $x\le y$, we use the notation
$$
A(x,y):=\exp\left(-\int_{x}^ya(\dd t)\right).
$$
\end{lemma}
\begin{theorem}[Schbath et al.]
\label{t.s}
For every $z$,
$$
r(z)=
\int_{x\le z\le y}\frac{J(x,z)\,J(z,y)}{J(x,y)}\,A(x,y)\,a(\dd x)\,a(\dd y).
$$
\end{theorem}

The contribution of the intensity measure $a$
in $r(z)$ corresponds to the product $A^-(z,\dd
x)\,A^+(z,\dd y)$.

A quick look at the ratio of the functions $J$ in the integral above
could lead to the erroneous conclusion that $r(z)$ is not well defined
when $J(x,y)$ is not always positive. (One knows that $J(x,y)$ is
positive when, for instance, the random variables $L_{t}$ are
uniformly integrable, and $c(\dd t)$ is uniformly bounded, that is,
when there exists a finite $\cho_{+}$ such that $c(\dd x)\le
\cho_{+}\,\dd x$.)  In fact, one can show that this ratio is at most
$1$ for any intensity $c(\dd t)$ and any distributions of the random
variables $L_{t}$, hence the formula for $r(z)$ in theorem~\ref{t.s}
is always valid.

We recall that, in the homogeneous case, the process of the clones has constant
intensity
$c(\dd x)=\cho\,\dd x$, the lengths of the clones are i.i.d.\ and 
distributed like a random variable $L$, and the
process of the anchors has constant intensity $a(\dd x)=\aho\,\dd x$.

\begin{corollary}[Arratia et al.]
\label{c.a}
In the homogeneous case with parameters $\cho$, $\aho$ and $L$,
$r(z)=\rhom$ does not depend on $z$
and its value is
$$
\rhom:=\int_{0}^{+\infty}\!\!\!\int_{0}^{+\infty}
\aho^2\,\ee^{-\aho (u+v)}\,\frac{J(u)\,J(v)}{J(u+v)}
\,\dd
u\,\dd v.
$$
Here, $J(u)$
is the probability that an interval of length $u$ is not covered by any
unique clone, hence
$$
J(u):=\exp\left(-\cho\,\int_u^{+\infty}\P(L\ge t)\,\dd t\right).
$$
\end{corollary}

When, furthermore, $L=\ell$ with full probability for a given positive real number 
$\ell$, 
Arratia et al.\ deduce from this the value of $\rhom$ as a function of $\ell$,
$\cho$
and $\aho$.

One gets the expression of $\rhom$ in corollary~\ref{c.a} from $r(z)$
in theorem~\ref{t.s}, using the change of variables $u=z-x$, $v=y-z$.

\section{Higher moments}
\label{s.hm}
Higher moments of the quantities introduced
above involve functionals of the processes that depend on 
more than one point. We first describe the computation of the
variance of the proportion of the real line which is occupied by the
ocean
in the general case,
then we consider the higher moments in the general case, and finally 
we prove precise asymptotics of the variance in the homogeneous case.

\subsection{Variance of the ocean proportion}

Recall that $r(z,z')$ is the probability that neither $z$ nor $z'$ are
covered by anchored clones.  Let $r_0(z,z')$, respectively
$r_1(z,z')$, respectively $r_2(z,z')$, denote the probability of the
same event, when the number of anchors between $z$ and $z'$ is $0$,
respectively $1$, respectively $2$ or more.  One can decompose each of
these events, according to the position of the first anchor to the
left of the interval $(z,z')$, which we call $x$ in the integrals
below, to the position of the first anchor to the right of $(z,z')$,
which we call $y$ in the integrals below, and to the positions of the
leftmost and rightmost anchors, if any, in the interval $(z,z')$,
which we call $s$ and $t$ in the integrals below.

Thus $r(z,z')=r_0(z,z')+r_1(z,z')+r_2(z,z')$ with, for $z\le z'$,
\begin{eqnarray*}
r_0(z,z')
& := & \int_{x\le z \le z'\le y}J(x \sa z,z'\sa y)\,B(\dd x,\dd y),
\\
r_1(z,z')
& := & \int_{x\le z\le s\le z'\le y}
J(x \sa z \sa s)\,J(s \sa z' \sa y)\,a(\dd s)\,B(\dd x,\dd y),
\\
r_2(z,z')
& := & \int_{x\le z\le s\le t\le z'\le y}
J(x \sa z \sa s)\,J(t \sa z'
\sa y)\,B(\dd x,\dd s)\,B(\dd t,\dd y).
\end{eqnarray*}
We mention that $r_{i}(z,z')$ is defined as an integral of dimension $i+2$, for $i=0$, $1$ or $2$.
We used the following notations. The two dimensional measure $B$
is defined on the subset $x\le y$ of $\R\times\R$
by the formula
$$
B(\dd x,\dd y):=A(x,y)\,a(\dd x)\,a(\dd y).
$$
For any $x\le z\le z'\le y$,
$$
J(x\sa z,z'\sa y):=\frac{J(x,z)\,J(z',y)}{J(x,y)},
\quad
J(x\sa z\sa y):=\frac{J(x,z)\,J(z,y)}{J(x,y)}.
$$ The quantities involved in the definitions above have the following
interpretations.  First, $\mathbf{1}\{x\le z\le y\}\,B(\dd x,\dd y)$
is the distribution of the couple formed by the positions of the
rightmost anchor to the left of $z$ and of the leftmost anchor to the
right of $z$. Second, $J(x\sa z\sa y)$ is the probability that $z$ is
not covered by an anchored clone when the closest anchor to the left
of $z$ is at $x$ and the closest anchor to the right of $z$ is at $y$.
Finally, $J(x\sa z,z'\sa y)$ is the probability that $z$ and $z'$ are
not covered by anchored clones when the closest anchor to the left of
$z$ is at $x$, the closest anchor to the right of $z'$ is at $y$, and
when there is no anchor between $z$ and $z'$.  Schbath's formula in
our theorem~\ref{t.s} reads
$$
r(z)=\int_{x\le z\le y}J(x \sa z\sa y)\,B(\dd x,\dd y).
$$
If one forgets the condition that $s\le t$ in the definition
of $r_2(z,z')$, one gets the product of the integrals over 
$(x,s)$ and over $(t,y)$, which are $r(z)$ and $r(z')$, respectively.
This implies our lemma~\ref{l.8} below. 

\begin{lemma}\label{l.8}
For any $z\le z'$,
$r_2(z,z')=r(z)\,r(z')-r_3(z,z')$ where the term
$r_3(z,z')$ is nonnegative and is
$$
r_3(z,z')
:= 
\int_{x\le z\le s,\,t\le z'\le y,\,s\ge t}
J(x \sa z \sa s)\,J(t \sa z'
\sa y)\,B(\dd x,\dd s)\,B(\dd t,\dd y).
$$
As a consequence, the variance $\sigma^2(O_G)$ of $O_G$ is
$$
\sigma^2(O_G)=\int_{0}^{G}\!\!\!\int_{0}^{G}(r_0+r_1-r_3)(z,z')\,\dd z\,\dd
z'.
$$
\end{lemma}

\subsection{Higher moments of the ocean proportion}

As mentioned above, one can adapt the technique used in the last section to study
the mean value of
any power of $O_G$. For
instance,
$$
\E(O_G^3)=\int_{0}^{G}\!\!\!\int_{0}^{G}\!\!\!\int_{0}^{G}r(z,z',z'')\,\dd z\,\dd z'\,\dd
z''.
$$
Thus, assuming for instance that $n=3$, one has to compute the $n$-point
function $r(z,z',z'')$.
First, one can assume by symmetry that $z\le z'\le z''$.
Let $x$ denote the position of the rightmost anchor to the left of
$z$, and $y$ the position of the leftmost anchor to the right of
$z''$. Let $s$ and $t$ denote the positions of the
leftmost and rightmost anchors in the interval $(z,z')$, and  $s'$ and
$t'$ 
the positions of the
leftmost and rightmost anchors in the interval $(z',z'')$, if these exist.

Then $r(z,z',z'')$ is $n!=6$ times the sum of $3^{n-1}=9$ 
terms $r_{i,i'}(z,z',z'')$.
Each term $r_{i,i'}(z,z',z'')$ corresponds to the number $i=0,1$ or
$2$ of anchors to be considered in the interval $(z,z')$ and to the
number
$i'=0,1$ or $2$ of anchors to be considered in the interval
$(z',z'')$.
Namely, no anchor at all, or a unique anchor, denoted by $s$ or by $s'$,
or two extremal anchors, denoted by $s$ and
$t$, or by $s'$ and $t'$. 

To take an example, consider the case $i=2$ and $i'=1$. This yields 
$r_{2,1}(z,z',z'')$ as the integral
$$
\int_{D_{2,1}} J(x \sa z \sa s)\,J(t \sa z'
\sa s')\,J(s'\sa z''\sa y)\,a(\dd s')\,B(\dd x,\dd s)\,B(\dd t,\dd y),
$$
where the domain of integration $D_{2,1}$ has dimension $5$ and is
defined by the inequalities
$$
x\le z\le s\le t\le z'\le s'\le
z''\le y.
$$
Likewise, if $i=0$ and $i'=2$, $r_{0,2}(z,z',z'')$ is the integral
$$
\int_{D_{0,2}} J(x \sa z, z'\sa s')\,J(t' \sa z''
\sa y)\,B(\dd x,\dd s')\,B(\dd t',\dd y),
$$
where the domain of integration $D_{0,2}$ has dimension $4$ and is defined by the inequalities
$$
x\le z,
\quad
z'\le s'\le t'\le
z''\le y.
$$ 
More generally, $\E(O_G^n)$ is the integral of
the $n$-point function $r(z_1,\ldots,z_n)$ on the domain $[0,G]^{n}$
with respect to the Lebesgue measure. 
For every $n$-tuple $z_1\le\cdots\le z_n$,
$r(z_1,\ldots,z_n)$ can be decomposed as a sum of $3^{n-1}$
contributions. Each of these contributions corresponds to the event
that each interval $[z_k,z_{k+1}]$ contains no anchor at all, or a
unique anchor, or at least two anchors.

\subsection{Variance in the homogeneous case}
\label{s.vhc}
In this section, we study the homogeneous case, when the intensity
measures are
$a(\dd t)=\aho\,\dd t$ and $c(\dd x)=\cho\,\dd x$, and the
distribution of the length $L_x$ of a clone does not depend on its
position $x$ and is the distribution of a random variable $L$. 
We recall that the distribution of the global process is
left invariant by the action of the translations. This 
implies that $r(z)=\rhom$ for every $z$, where the value of $\rhom$
is given in corollary~\ref{c.a}. Hence, 
$$\E(O_G)=G\,\rhom.
$$
Since $(z,z')\mapsto r(z,z')-r(z)\,r(z')$ is a symmetric function,
$\sigma^2(O_G)$ is twice an integral over $z'\ge z$. 
Likewise, the invariance by the translations implies that 
$r(z,z')=r(0,z'-z)$ for every $z$ and $z'$.
Introducing $\bar r_i(z):=r_i(0,z)$, one is left with twice some
integrals of the functions $\bar r_i(z)$ over
$z$ in $[0,G]$, namely
$$
\sigma^2(O_G)=2\int_{0}^{G}(G-z)\,
(\bar r_0(z)+\bar r_1(z)-\bar r_3(z))\,\dd z.
$$
The values of the quantities $\bar r_i(z)$ for every nonnegative $z$ are 
\begin{eqnarray*}
\bar r_0(z)
& = & \int_{x, y\ge
  0}\aho^2\,\ee^{-\aho(x+y+z)}\,\frac{J(x)\,J(y)}{J(x+y+z)}\,\dd x\,\dd y,
\\
\bar r_1(z)
& = & \int_{x, y\ge
  0,\, 0\le t\le
  z}\aho^3\,\ee^{-\aho(x+y+z)}\,\frac{J(x)\,J(t)\,J(z-t)\,J(y)}{J(x+t)\,J(z-t+y)}
\,\dd x\,\dd y\,\dd t,
\\
\bar r_3(z)
& = & \int_{x,y,s, t\ge 0, s+t\ge
  z}\aho^4\,\ee^{-\aho(x+y+s+t)}\,\frac{J(x)\,J(t)\,J(s)\,J(y)}{J(x+t)\,J(s+y)}
\,\dd x\,\dd y\,\dd s\,\dd t.
\end{eqnarray*}
We mention that $\bar r_{0}(z)$, respectively $\bar r_{1}(z)$, respectively $\bar r_{3}(z)$, is defined as 
an integral of dimension $2$, respectively $3$, respectively $4$.

Using the fact that the function $x\mapsto J(x)$ is nondecreasing, one can bound each $\bar r_{i}(z)$ as follows:
\begin{eqnarray*}
\bar r_0(z) &\le & \ee^{-\aho\,z}\,j(\aho),
\\
\bar r_1(z) & \le & \aho\,z\,\ee^{-\aho\,z}\,j(\aho)^2,
\\
\bar r_3(z) & \le & (1+\aho\,z)\,\ee^{-\aho\,z}\,j(\aho)^2,
\end{eqnarray*}
with the notation
$$
j(\aho):=\int_{0}^{+\infty}\aho\,\ee^{-\aho\,x}\,J(x)\,\dd x.
$$ To prove the upper bound of $\bar r_{0}(z)$, one uses the fact that
$J(y)\le J(x+y+z)$, and one performs the integration of the upper
bound.  Likewise, to prove the upper bound of $\bar r_{1}(z)$, one
uses the facts that $J(t)\le J(x+t)$ and $J(z-t)\le J(z-t+y)$, and one
performs the integration of the upper bound.  Finally, to prove the
upper bound of $\bar r_{3}(z)$, one uses the facts that $J(t)\le
J(x+t)$ and $J(s)\le J(s+y)$, and one performs the integration of the
upper bound. In this last case, this yields
$$
\bar r_{3}(z)
\le
j(\aho)^{2}\int_{s, t\ge 0, s+t\ge
  z}\aho^2\,\ee^{-\aho(s+t)}\,\dd s\,\dd t,
$$
and the last double integral is indeed $(1+\aho\,z)\,\ee^{-\aho\,z}$.

Since $J(x)\le1$, $j(\aho)\le1$.
Furthermore, the limit of $J$ at infinity is $1$, hence
$\bar r_0(z)\sim\ee^{-\aho\,z}\,j(\aho)^2$ at infinity.
Let
$$
\sigma^2_i(G):=\int_{0}^G2(G-z)\,\bar r_i(z)\,\dd z.
$$
From the bounds on the three functions $\bar r_i$ which are stated above,
it is not difficult to prove that, when $G\to\infty$,
$$
\sigma^2_i(G)=\nu_i\,G-\lambda_i+\tau_i(G),
$$ 
where $\tau_i(G)=o(1)$ for $i=0,1$ and $3$. 
More specifically, these bounds imply that the
numbers $\nu_i$ and $\lambda_i$, defined as 
$$
\nu_i:=\int_{0}^{+\infty}2\bar r_i(z)\,\dd z,
\qquad
\lambda_i:=\int_{0}^{+\infty}2z\,\bar r_i(z)\,\dd z.
$$
are indeed finite and positive, and simple computations show that 
$$
\tau_i(G):=\int_{G}^{+\infty}2(z-G)\,\bar r_i(z)\,\dd z.
$$
Introduce
$\tau(G):=\tau_0(G)+\tau_1(G)-\tau_3(G)$.
Since each $\tau_i(G)$ is nonnegative, $|\tau(G)|$ is at most the
maximum of $\tau_0(G)+\tau_1(G)$ and $\tau_3(G)$. Since
$j(\alpha)\le1$, our bounds on the three functions $\bar r_i$
imply that
$$
|\tau(G)|\le\int_G^{+\infty}2(z-G)\,(1+\aho\,z)\,\ee^{-\aho
  z}\,\dd z.
$$
Performing the integration, one gets
$$
|\tau(G)|\le2\aho^{-2}\,(3+\aho\,G)\,\ee^{-\aho G}.
$$
Finally, when $G\to\infty$,
$\tau(G)=O(G\,\ee^{-\aho\,G})$.

Assume now that $L\le\ell$ almost surely, for a
finite $\ell$.
This means that the intensity measure of the
global Poisson process on $\R\times\R^+$ puts no
mass on the set $\R\times(\ell,+\infty)$.
Assume that $z$ and $z'$ are such that $|z-z'|>\ell$. Then $I(z)\cap
I(z')$ contains only clones $(x,t)$ such that both points $z$ and $z'$
belong to $[x-t,x]$, hence in particular, such that $t>\ell$.
Since $I(z)\cap
I(z')$ is a subset of $\R\times(\ell,+\infty)$, its intensity measure
must be zero. Thus, the events $\{z\in
O\}$ and $\{z'\in O\}$ are in fact measurable with respect to 
the truncated cones of influence $I(z)\cap(\R\times[0,\ell])$ and
$I(z')\cap(\R\times[0,\ell])$, respectively. Since these two subsets of
$\R\times\R^+$ are disjoint, $\{z\in
O\}$ and $\{z'\in O\}$ are independent events. 

Finally, if $L\le\ell$ almost surely, $r(z,z')=r(z)\,r(z')$
as soon as $z$ and $z'$ are such that $|z-z'|>\ell$, 
hence $\bar r_0(z)+\bar r_1(z)-\bar r_3(z)=0$ for every $z>\ell$, and
$\tau(G)=0$ for every $G\ge\ell$.

Proposition~\ref{p.ss} below summarizes the results of this section.
\begin{proposition}
\label{p.ss}
{\bf (1)}
Let $\nu:=\nu_0+\nu_1-\nu_3$ and
$\lambda:=\lambda_0+\lambda_1-\lambda_3$.
When $G\to\infty$,
$$
\sigma^2(O_G)=\nu\,G-\lambda+o(1).
$$ 
{\bf (2)}
Assume that $L\le\ell$ almost surely for a
finite $\ell$. Then, for every 
$G\ge\ell$,
$$
\sigma^2(O_G)=\nu\,G-\lambda.
$$
\end{proposition}

\section{Functional invariance in the homogeneous case}
\label{s.ihc}
Our main task in this section is to prove that $\nu$ is positive, that
is, not zero.
We do this, first, in the limit $\cho\to0$ of a vanishing number of
clones, 
then in the
general case.
Our techniques also yield upper and lower bounds of the mean value and of the
variance of $O_G$ when the intensities are not constant.
Finally, we prove the functional invariance result of theorem~\ref{t.a}.

\subsection{Variance for vanishing clones}

\begin{proposition}[Homogeneous case]
\label{p.vvc}
Fix the distribution of 
  $L$ and the value of $\aho$. Then, if $\cho$ is small enough, 
  $\nu$ is positive. More precisely, when $\cho\to0$,
$$
\nu=\aho^{-2}\,\E(\varphi(\aho\,L))\,\cho+o(\cho),
$$ 
where the function
  $x\mapsto\varphi(x)$ is explicit, positive on 
$x>0$, and given by the formula
$$
\varphi(x):=x-1+\ee^{-x}\,(1-x^2/2).
$$
\end{proposition}
\proof{Proof}
If $\cho=0$, $j(\aho)=1$ and
$\bar r_i(z)=r^*_i(z)$, with
$$
r^*_0(z):=\ee^{-\aho\,z},
\quad
r^*_1(z):=\aho\,z\,\ee^{-\aho\,z},
\quad
r^*_3(z):=(1+\aho\,z)\,\ee^{-\aho\,z},
$$ 
hence $r^*_0+r^*_1-r^*_3$ is identically zero. (Besides, when $\cho=0$,
$O_G$ is almost surely zero.)
We now show that the first derivative of $\nu$ with respect to $\cho$ at $\cho=0^+$
is positive.

When $\cho=o(1)$, $J(x)=1-\cho\,H(x)+o(\cho)$ with
$$
H(x):=\displaystyle\int_{x}^{+\infty}\P(L\ge
t)\,\dd t.
$$
This implies that
$\bar r_i(z)=r^*_i(z)+\cho\,s_i(z)+o(\cho)$, for some explicit functions $s_i(z)$.
Introducing $w_i:=\displaystyle\int_{0}^{+\infty}s_i(z)\,\dd z$ and 
$w:=w_0+w_1-w_3$, one gets
$\nu=\cho\,w+o(\cho)$. 
For instance,
$$
w_0=\int_{x,y,z\ge0}\aho^2\,\ee^{-\aho\,(x+y+z)}\,\{H(x+y+z)-
H(x)-H(y)\}
\,\dd
  x\,\dd y\,\dd z,
$$
and similar expressions of $w_1$ 
and $w_3$ obtain.
After some tedious but simple computations, one gets
$$
w_0=h_2-2h_0,
\quad
w_1=2h_1-4h_0,
\quad
w_3=2h_2-6h_0,
$$ 
where, for every nonnegative integer $n$, the value of $h_{n}$ is given by
$$
h_n:=\displaystyle\int_{0}^{+\infty}\frac{(\aho\,x)^n}{n!}\,\ee^{-\aho\,x}\,H(x)\,\dd x.
$$
Summing up these three contributions yields $w=2h_1-h_2$. Converting everything back in terms of the distribution
of $L$, one finally gets
$$
w=\aho^{-2}\,\E(\varphi(\aho\,L)),
$$
where $\varphi$ is given in the statement of the proposition above.
It happens that $\psi(x):=\ee^x\,\varphi(x)$ defines a function $\psi$
such that $\psi(0)=0$ and whose derivative $\psi'(x)=x\,(\ee^x-1)$ is
obviously positive for every positive $x$.
Thus $\varphi(x)$  is
positive for every positive $x$, and $w$ is positive for every
distribution of $L$, except in the degenerate case when
$L=0$ almost surely. This proves that $\nu$ is positive for
small values of $\cho$.
\endproof

\begin{rbox}
Other limiting cases are possible. Recall that $\E(O_G)=\rhom\,G$
for every nonnegative $G$, and that
$\sigma^2(O_G)\sim\nu\,G$ when $G\to\infty$.
\begin{enumerate}
\item
If $\E(L^3)=o(1)$, then 
$\nu\sim\frac13\,\aho\,\cho\,\E(L^3)$.
\item
If $\cho=o(1)$, then $(1-\rhom)\sim \cho\,\E(L\,\ee^{-\aho L})$.
\item
If
$\E(L)=o(1)$, then $(1-\rhom)\sim \cho\,\E(L)$. 
\end{enumerate}
Note that this last result
does not depend on the value of $\aho$.
\end{rbox}

\subsection{Positive dependence}
\label{ss.pd}
Proposition \ref{p.3} below deals with possibly inhomogeneous
processes.

\begin{proposition}[General case]
\label{p.3}
For every Borel sets $Z$ and $Z'$,
$$
\P(Z\cup Z'\subset O)\ge\P(Z\subset O)\,\P(Z'\subset O).
$$
In particular, $r(z,z')\ge r(z)\,r(z')$ for every $z$ and $z'$.
\end{proposition}
Corollary \ref{c.3} is a direct consequence of this proposition and of
the expression of $\sigma^2(O_G)$ in section \ref{s.vhc}.
\begin{corollary}[Homogeneous case]
\label{c.3}
For every
nonzero intensities $\cho$ and $\aho$ and every nonzero $L$, 
the constants $\nu$ and $\lambda$
 are positive and the function $G\mapsto\tau(G)$ is nonnegative.
In particular, for every $G$,
$$
\nu\,G-\lambda\le\sigma^2(O_G)\le\nu\,G.
$$
Hence $\sigma^2(O_G)\sim\nu\,G$ when $G\to\infty$.
Furthermore, the following properties hold.
The function $G\mapsto\sigma^2(O_G)$ is increasing and
convex. When $G\to0$,  $\sigma^2(O_G)\sim\rhom\,(1-\rhom)\,G^2$.
When $G\to\infty$, $\sigma^2(O_G)=\nu\,G-\lambda+o(1)$.
\end{corollary}
\proof{Proof of corollary \ref{c.3}}
As regards $\nu$, recall from section~\ref{s.vhc} that,
in the homogeneous case,
$$
\nu=\int_0^{+\infty}2(r(0,z)-\rhom^2)\,\dd z.
$$
Since $0<\rhom<1$, $r(0,0)=r(0)=\rhom>\rhom^2$. 
Furthermore, one can deduce from section~\ref{s.hm} an expression of
$r(0,z)$ 
from the
formulas which give $r_i(z,z')$ for $i=0$, $1$ and $2$. The integrals
involved are continuous with respect to $z$ and $z'$ because the
functions $J$ involved in these integrals are, 
and because obvious domination properties hold. Finally,
$r(0,z)>\rhom^2$ for every nonnegative $z$ in a neighborhood of $0$,
and  $r(0,z)\ge\rhom^2$ for every nonnegative $z$. This implies that
$\nu>0$.

 The proofs that $\lambda$ is positive and that $\tau(G)$ is nonnegative
 are similar.

The equivalent of $\sigma^2(O_G)$ when $G\to0$ stems from  the fact
that $r(0,z)\to\rhom$ when $z\to0$ and from the exact
formula
$$
\sigma^2(O_G)=\int_0^G2(G-z)\,(r(0,z)-\rhom^2)\,\dd z.
$$
Finally, this formula and the fact that $r(0,z)\ge\rhom^2$
also yield the fact that 
the function $G\mapsto\sigma^2(O_G)$ is increasing and
convex, since the derivative of this function is
$$
\int_0^G2(r(0,z)-\rhom^2)\,\dd z.
$$ 
\endproof
\proof{Proof of proposition \ref{p.3}}  
For any Borel set $Z$, 
$\{Z\subset O\}$ is a nonincreasing event, with respect to the global Poisson process
introduced in section~\ref{ss.ca}. To see this, note that, if one adds some anchors and/or some clones to
a given configuration, the union $\R\setminus O$ of the anchored
islands does not decrease hence the indicator function of the event 
$\{Z\subset O\}$ does not increase.
Thus, our proposition is a direct consequence of the
Fortuin-Kasteleyn-Ginibre (FKG) inequality
$$
\P(A\cap A')\ge\P(A)\,\P(A'),
$$ 
applied to the nonincreasing events $A:=\{Z\subset O\}$ and 
$A':=\{Z'\subset O\}$, see Roy~(1991) for instance.
\endproof

\subsection{Bounds in the general case}

In the inhomogeneous case, minimal assumptions on $c(\dd x)$ and
$a(\dd x)$ yield upper and lower bounds on $\E(O_G)$ and
$\sigma^2(O_G)$, as we now show.  In this section, we assume that the
intensities of the processes of the clones and of the anchors are
uniformly bounded. Hence, $a(\dd x)$ and $c(\dd x)$ are absolutely
continuous with respect to the Lebesgue measure and there exists
finite positive constants $\aho_{\pm}$ and $\cho_{\pm}$ such that
$$
\aho_-\,\dd x\le a(\dd x)\le \aho_+\,\dd x,
\quad
\cho_-\,\dd x\le c(\dd x)\le \cho_+\,\dd x.
$$
We assume furthermore that the lengths $L_x$ of the clones are uniformly
stochastically bounded from above and from below.
This means that there exists nonnegative random variables $L_{\pm}$
such that $L_+$ is integrable, such that $L_-$ is not almost surely zero, and,
such that, for every $x$ and $t$,
$$
\P(L_-\ge t)\le\P(L_x\ge t)\le\P(L_+\ge t).
$$
In particular, the family $(L_x)_x$ must be uniformly integrable.

\begin{proposition}
\label{p.unif}
The assumptions above imply that there exists positive
constants $\rhom_{\pm}<1$ and finite positive constants $\nu_{\pm}$ such that, for every $G$,
$$
\rhom_-\,G\le\E(O_G)\le\rhom_+\,G,
\quad
\nu_-\,G\le\sigma^2(O_G)\le\nu_+\,G.
$$
\end{proposition}

In these inequalities, 
$\rhom_-$ corresponds to the homogeneous case of parameters $\cho_+$,
$\aho_+$ and 
$L_+$, and $\rhom_+$ 
 to the homogeneous case of parameters $\cho_-$, $\aho_-$ and $L_-$.
As regards the variance, the dependence is not so straightforward, at
least the dependance that our techniques yield.
The parameter $\nu_+$ that we exhibit depends on $\aho_-$ alone, a result
which may seem surprising, and the parameter
$\nu_-$ depends on $\cho_+$, $\rhom_+$, and $\rhom_-$. 

\proof{Proof}
The bounds on $\E(O_G)$ would follow from the fact that
$$
\rhom_-\le r(z)\le\rhom_+,
$$ for any $z$ and for positive $\rhom_{\pm}<1$.  Such bounds on
$r(z)$ themselves stem from the fact that the distribution of the
ocean, as a random subset of the real line, is nonincreasing with
respect to the intensities of the processes of the clones and of the
anchors.  Hence, by a coupling argument, the value of $\E(O_G)$ lies
between its value for the homogeneous processes of densities $\aho_+$
and $\cho_+$ on the one hand, and $\aho_-$ and $\cho_-$ on the other
hand, the distributions of the lengths $L_x$ being fixed.

We now examine the
influence of the distributions of the lengths.  Once again by a
coupling argument, the uniform
replacement of the distributions of the lengths $L_x$ by the
distribution of $L_+$ yields longer clones, hence longer islands,
hence a stochastically smaller ocean. This proves the lower bound of
$\E(O_G)$.  Comparison with $L_-$ yields the upper bound.

Our proof of the lower bound of $\sigma^2(O_G)$ goes as follows.
One knows that
$$
\sigma^2(O_G)=\int_{0}^{G}\!\!\!\int_{0}^{G}(r(z,z')-r(z)\,r(z'))\,\dd z\,\dd
z',
$$
and that the expression $r(z,z')-r(z)\,r(z')$ is nonnegative for every
$z$ and $z'$.
Assume that there exists positive $\delta$ and $\varepsilon$ such that, for every
$z$ and $z'$ such that $|z-z'|\le\varepsilon$,
$$
r(z,z')-r(z)\,r(z')\ge\delta.
$$
The lower bound of $\sigma^2(O_G)$ would follow.
Now, for every $z\le z'$,
if $z'$ is in $O$ and if there is no right end of clone in $[z,z']$,
then $z$ is in $O$.
Hence,
$$
r(z')=r(z,z')+\P(z\notin O,\,z'\in O)\le r(z,z')+\P(B),
$$
with $B:=\{\ccp\cap([z,z']\times\R^+)\neq\vide\}$.
By definition of the intensity of the Poisson process $\ccp$,
$$
\P(B)=1-\ee^{-c([z,z'])}\le c([z,z'])\le \cho_+\,(z'-z).
$$
Since $r(z')\ge\rhom_-$ and $r(z)\le\rhom_+$, this proves the lower bound
$$
r(z,z')-r(z)\,r(z')\ge(1-\rhom_+)\,\rhom_--\cho_+\,(z'-z).
$$
This in turn shows the desired inequality for $z\le z'$ and $z'-z$
small enough.

As regards the upper bound, it is enough to bound from above the
integrals of $r_0(z,z')$ and $r_1(z,z')$, since $r_3(z,z')$ is nonnegative.
In the expression of $r_0(z,z')$, for every fixed values of $x$ and $y$,
$J(x\sa z,z'\sa y)$ is a nonincreasing function of the distributions
of the lengths $L_x$ and of the intensity of the
clones, since having
more clones and longer clones only makes the ocean smaller.
Thus  $r_0(z,z')$ is bounded from above by its value when one replaces
$c(\dd t)$ by
$\cho_-\,\dd t$ and the distribution of every $L_x$ by the distribution of
$L_-$.
Likewise, the interpretation of $B(\dd x,\dd y)$ as the joint
distribution of the positions of the rightmost anchor to the left of
$z$ and of the leftmost anchor to the right of $z$, and a coupling
between two processes of anchors with comparable intensities, show that
the anchors become stochastically more distant from $z$ when one
replaces 
$a(\dd t)$ by the smaller intensity $\aho_-\,\dd t$. Hence the
probability that $z$ is not covered by an anchored clone cannot decrease.
Thus, replacing  $a(\dd t)$  by
$\aho_-\,\dd t$ cannot make $r_0(z,z')$ decrease. 

Finally, the
contribution of $r_0$ in the value of 
$\sigma^2(O_G)$ is bounded from above by its value in the homogeneous
case which uses the values $\aho_-$, $\cho_-$ and $L_-$,
that is, for instance, by $2G/\aho_-$.
Likewise, the
contribution of $r_1$ to the value of 
$\sigma^2(O_G)$ is at most $2G/\aho_-$.
This yields the desired upper bound with  $\nu_+:=4/\aho_-$.
\endproof

\begin{rbox}
Alternatively, when $L_z\le\ell$ almost surely and for every $z$,
recall from the end of section~\ref{s.vhc} that 
$r(z,z')=r(z)\,r(z')$ as soon as $|z-z'|>\ell$, 
hence $\sigma^2(O_G)$ is at most the area of the part of the
square $[0,G]^2$ inside the diagonal strip
$|z-z'|\le\ell$, that is, at most $2\ell G-\ell^2$ when
$G\ge\ell$, and $\sigma^2(O_G)$ is at most $G^2$ for every $G$.
Hence $\sigma^2(O_G)\le2\ell G$ for every $G$.
\end{rbox}
\begin{rbox}
One can adapt the proofs in this section to some cases when the intensities
of the clones and of the anchors are zero in some places, as long as the
intensities stay bounded from below on regions which are spread out enough.
\end{rbox}

\subsection{Convergence in distribution}

We first explain how one could prove the convergence of the moments by
elementary techniques, then we show that general invariance results apply,
which yield directly the desired convergence.

\subsubsection{Method of moments}

Assume first that $L\le\ell$ almost surely. Then, a crucial
remark from the end of section~\ref{s.vhc} 
is that the events $\{Z\subset O\}$ and
$\{Z'\subset O\}$ are independent as soon as the distance between every point in $Z$ and
every point in $Z'$ is at least $\ell$.
Furthermore,
$$
\E((O_G-\rhom\,G)^n)=\int_{Z\in[0,G]^n}\pi(Z)\,\dd Z,
\quad
\pi(Z):=\prod_{z\in Z}(\un\{z\in O\}-\rhom).
$$
If $Z=Z'\cup Z''$ with $|z'-z''|\ge\ell$ for every
$z'\in Z'$ and every $z''\in Z''$, one gets
$\E(\pi(Z))=\E(\pi(Z'))\,\E(\pi(Z''))$.

For instance, if $n=3$, 
every nontrivial partition of $Z$ includes at least one singleton
hence $\E(\pi(Z))$ is zero except when all the distances between
the nonempty subsets
of $Z$ are at most $\ell$. Ordering the points $z$, $z'$ and $z''$, we are left with the domain 
$$z\le z'\le
z+\ell,\quad
z'\le z''\le z'+\ell,
$$
whose volume is at most $\ell^2G$. Hence $\E((O_G-\rhom\,G)^3)=O(G)$.

If $n=4$, the only difference with the $n=3$ case is due to
the partitions of $Z$ into two pairs $Z'$ and $Z''$. These contribute to
the result even when  the distance from $Z'$ to $Z''$ is large. Every
such 
$\E(\pi(Z'))$ and $\E(\pi(Z''))$
is $O(G)$, hence $\E((O_G-\rhom\,G)^4)$ is $O(G^2)$.

Likewise, for every positive integer $k$,  the moments $\E((O_G-\rhom\,G)^{2k})$ and $\E((O_G-\rhom\,G)^{2k+1})$ are
both $O(G^k)$.

One can also compute the asymptotics of the moments of $O_G$ as
$G\to\infty$.  To do this, one starts from the expression of
$\E((O_G-\rhom\,G)^{2n})$ as the integral of $\E(\pi(Z))$ over the
points $Z$ in $[0,G]^{2n}$. When there exists a partition of $Z$ into
two parts $Z'$ and $Z''$ at a distance at least $\ell$, $\E(\pi(Z))$
is the product $\E(\pi(Z'))\,\E(\pi(Z''))$.  The remaining points $Z$
span a volume in $[0,G]^{2n}$ which is $o(G^n)$, hence they contribute
to a vanishing part of the asymptotics.

This yields recursions between the asymptotic moment of degree $2n$
and 
the asymptotic moments of even degrees at most $2n-2$.
One can deduce from these recursions the convergence of the moments of 
$(O_G-\rhom\,G)/\sqrt{G}$ to the moments of a Gaussian random variable.

Finally, one could adapt this strategy to the case where $L$ is
unbounded, thus reaching the same conclusion.

\subsection{Direct method}

A stronger conclusion obtains directly from classical results by Doukhan
et al.\ (1994), for every square integrable  $L$.
To see this, introduce for every integer $n$, the random variable
$$X_n:=O([n,n+1])-\rhom.$$
Let $\cf_n$ denote the $\sigma$-algebra 
generated by the collection $(X_i)_{i\le n}$, 
and let $\cg_n$ denote the $\sigma$-algebra generated by the collection $(X_i)_{i\ge n}$.
The sequence $(X_n)_{n}$ is generated by the action of the shift
$$
\vartheta:(x,t)\mapsto (x+1,t),
$$ 
on $\R\times M$, since $X_n=X_0\circ\vartheta^n$ for every integer $n$.
The strong mixing coefficients $\alpha_n$ associated to the stationary sequence
$(X_n)_{n}$ are defined, for any integer $n\ge0$, by
$$
\alpha_n:=\sup\{\P(B\cap
B')-\P(B)\,\P(B')\,;\,B\in\cf_0,\,B'\in\cg_n\}.
$$
Since $|X_0|\le1$ almost surely, the condition in Doukhan
et al.\ (1994) reduces to the summability of the series of general
term $\alpha_n$.
Neglecting the influence of the anchors does not decrease 
the value of $\alpha_n$. Thus $\alpha_n\le\P(B_n)$, where
$B_n$ is the event that at least one clone covers both points
$0$ and 
$n$.
One can bound each $\P(B_{n})$ as follows:
$$
\P(B_n)=1-J(n)\le\int_{n}^{+\infty}\cho\,\P(L\ge t)\,\dd t.
$$ This shows that the sequence of general term $\P(B_n)$ is summable
as soon as $L$ is square integrable. (In fact, this sequence is
summable if and only if $L$ is square integrable, proof omitted.)
This shows that the functional invariance stated in theorem \ref{t.a}
holds, at least for the processes $\mao_G$ such that $G$ is an
integer. The general case is an easy consequence, since $O_G$ depends
on $G$ in a monotone way.

Equivalently, one can write directly $O_G$ as
$$O_G=\rhom\,G+\int_0^GY_x\,\dd x,$$ where the stationary centered family
$Y_x:=\un\{x\in O\}-\rhom$ is indexed by the real numbers $x$. The same
conclusion obtains.

\bigskip

Universit\'e Claude Bernard Lyon 1 \\ 
Institut Camille Jordan UMR 5208 
\\
Domaine de Gerland \\
50, avenue Tony-Garnier \\ 69366 Lyon Cedex 07 (France) 

{\tt
Didier.Piau@univ-lyon1.fr} \\ {\tt http://lapcs.univ-lyon1.fr}


\begin{thebibliography}{1}

\bibitem{altw}
Richard Arratia, Eric S.\ Lander, Simon Tavar\'e, and Michael S.\ Waterman.
Genomic mapping by anchoring random clones: a mathematical analysis. 
{\em Genomics\/} {\bf 11}, 806-827 (1991). 

\bibitem{dmr}
Paul Doukhan, Pascal Massart, and Emmanuel Rio.
The functional central limit theorem for strongly mixing
processes.
{\em Annales de l'Institut Henri Poincar\'e
Probabilit\'es et Statistique\/} {\bf 30}, 63-82 (1994).

\bibitem{ebddme}
Warren J.\ Ewens, Callum J.\ Bell, Peter J.\ Donnelly, Patrick Dunn, 
Emilia Matallana, and Joseph R.\ Ecker.
Genome mapping with anchored clones: theoretical aspects.
{\em Genomics\/} {\bf 11}, 799-805 (1991). 

\bibitem{gri}
Andrei V.\ Grigoriev. Theoretical predictions and experimental
observations of genomic mapping  
by anchoring random clones. 
{\em Genomics\/} {\bf 15}, 311-316 (1993). 

\bibitem{lw}
Eric S.\ Lander, Michael S.\ Waterman. Genomic mapping by fingerprinting random clones: a mathematical analysis. 
{\em Genomics\/} {\bf 2}, 231-239 (1988).

\bibitem{r}
Rahul Roy.
Percolation of Poisson sticks on the plane. 
{\em Probability Theory and Related Fields\/} {\bf 89}, 503-517 (1991).
 
\bibitem{s}
Sophie Schbath. Coverage processes in physical mapping by anchoring random
clones. {\em Journal of Computational Biology\/} {\bf 4}, 61-82 (1997). 
 
\bibitem{sbt}
Sophie Schbath, Nathalie Bossard, and Simon Tavar\'e. The effect of non-homogeneous
clone length distribution on the progress of an STS mapping
project. {\em Journal of Computational Biology\/} {\bf 7}, 47-58 (2000). 

\end{thebibliography}
\end{document}